\documentclass[12pt,reqno,a4wide]{amsart}

\usepackage{tikz} 
\usepackage{calrsfs}
\usepackage{mathrsfs}
\usetikzlibrary{decorations.pathmorphing}
\usetikzlibrary{arrows.meta}
\usepgflibrary {shadings}
\usepackage{array}
\usepackage{xcolor}

\usepackage{comment}
\usepackage{hyperref}

\usetikzlibrary{decorations.pathreplacing}
\usetikzlibrary{fit}

\usepackage{pgfplots}

\allowdisplaybreaks

\catcode`,\active

\catcode`\,12

\begin{document}
\bibliographystyle{plain}

%
%

	\title[Abdelkader's random walk model]
	{Generating functions and Abdelkader's random walk model}

	\author[H. Prodinger ]{Helmut Prodinger }
	\address{Department of Mathematics, University of Stellenbosch 7602, Stellenbosch, South Africa
	and
NITheCS (National Institute for
Theoretical and Computational Sciences), South Africa.}
	\email{hproding@sun.ac.za}

	\keywords{Random walks with resets, Moran walks,   generating functions }
\subjclass{05A15}

	\begin{abstract}
We link questions by Abdelkader about a class of random walks to \emph{Moran walks}.

	\end{abstract}
	
	\subjclass[2020]{05A15}

\maketitle

\section{Introduction}

M. Abdelkader studied a random walk model where one can advance by one step with
probability $p$ or fall down to the $x$-axis with probability $q=1-p$. The height of an
$n$ step random walk of this type is the the maximal ordinate achieved during the walk.

He asked:

\begin{itemize}
	\item [1.] Can we find a closed form of the probability-generating function of the height?\\
	\item [2.] Can we explicitly calculate the mean and variance of the height statistics using the
	probability-generating function of $H_n$?
\end{itemize}

  \section{The Moran model}
  
  Abdelkader's model is known as the \emph{Moran model} and was recently given a very thorough study in \cite{Moran}.

Consequently we refine ourselves to a few comments, by shedding some light on the machinery of generating functions and suggest to read the long cited paper.

Let $q=1-p$. The admissible words can be decomposed according to their return to the $x$-axis, and a sequence of
 possible up-steps at the end.
One sojourn is enumerated by
\begin{equation*}
X=\sum_{n\ge1}p^nz^{n+1}q=\frac{pqz^2}{1-pz};
\end{equation*}
an arbitrary sequence of sojourns is enumerated by
\begin{equation*}
\frac1{1-X}=\frac{1-pz}{1-pz-pqz^2}
\end{equation*}
and with an possible sequence of up-steps at the end
\begin{equation*}
\frac{1-pz}{1-pz-pqz^2}\sum_{n\ge0}p^nz^n=\frac1{1-pz-pqz^2}.
\end{equation*}
We need an asymptotic equivalent of the coefficient of $z^n$ in this, written as $[z^n]\frac1{1-pz-pqz^2}$.
Since we just have a rational function, we even have an explicit form (aka `Binet' form):
\begin{equation*}
[z^n]\frac1{1-pz-pqz^2}=\frac1{\sqrt{p(p+4q)}}\biggl[\biggl(\frac{p+\sqrt{p(p+4q)}}{2}\biggr)^{n+1}-\biggl(\frac{p-\sqrt{p(p+4q)}}{2}\biggr)^{n+1}\biggr].
\end{equation*}
In terms of asymptotics, only the first term is relevant, and we have
\begin{equation*}
	[z^n]\frac1{1-pz-pqz^2}\sim\frac1{\sqrt{p(p+4q)}}\biggl(\frac{p+\sqrt{p(p+4q)}}{2}\biggr)^{n+1}.
\end{equation*}

Now we move to generating functions with bounded height:  any such path must have all ordinates $\le H$. Then we get
\begin{equation*}
	X=\sum_{n=1}^Hp^nz^{n+1}q=\frac{pqz^2}{1-pz}-\frac{qp^{H+1}z^{H+2}}{1-pz},
\end{equation*}
\begin{equation*}
\frac1{1-X}=\frac{1-pz}{1-pz-pqz^2+qp^{H+1}z^{H+2}},
\end{equation*}
\begin{equation*}
F^{[\le H]}:=	\frac1{1-X}\sum_{n=0}^Hp^nz^n=\frac{1-p^{H+1}z^{H+1}}{1-pz-pqz^2+qp^{H+1}z^{H+2}}.
\end{equation*}
An explicit form is not obvious, so we switch to asymptotics, and we consider the smallest zero of the denominator, which we call $\varepsilon=\varepsilon_H$.
We follow the classic bootstrapping method as described in \cite{Knuth78},
\begin{equation*}
1-p\varepsilon-pq\varepsilon^2\sim0;
\end{equation*}
so we get a first approximation
\begin{equation*}
\varepsilon\sim\frac{p+\sqrt{p(p+4q)}}{2pq}.
\end{equation*}
For the next approximation, we include
\begin{equation*}
qp^{H+1}\varepsilon^{H+2}\sim \frac{q}{p}\biggl(\frac{p+\sqrt{p(p+4q)}}{2q}\biggr)^{H+2}.
\end{equation*}
The next equation is thus
\begin{equation*}
	1-p\varepsilon-pq\varepsilon^2+ \frac{q}{p}\biggl(\frac{p+\sqrt{p(p+4q)}}{2q}\biggr)^{H+2}\sim0.
\end{equation*}
Ignoring irrelevant terms, we find a better expansion:
\begin{equation*}
\overline{\varepsilon}=\varepsilon+\frac{q}{p\sqrt{p(p+4q)}}(p\varepsilon)^{H+2}
\end{equation*}
\begin{equation*}
\frac{1-p^{H+1}z^{H+1}}{1-pz-pqz^2+qp^{H+1}z^{H+2}}\sim\frac{1}{-p-2{\varepsilon}pq}\frac1{z-\overline{\varepsilon}}
\sim\frac{1}{p(1+2{\varepsilon}q)\varepsilon}\frac1{1-z/\overline{\varepsilon}}
\end{equation*}
and
\begin{equation*}
[z^n]\frac{1-p^{H+1}z^{H+1}}{1-pz-pqz^2+qp^{H+1}z^{H+2}}\sim\frac{1}{p(1+2{\varepsilon}q)}\frac1{\varepsilon^{n+1}}
\Bigl(1-\frac{q}{\sqrt{p(p+4q)}}(p\varepsilon)^{H+1}\,\Bigr)^n
\end{equation*}
Of course, $H\to\infty$ means no restrictions, and we get then
\begin{equation*}
	[z^n]\frac{1}{1-pz-pqz^2}\sim\frac{1}{p(1+2{\varepsilon}q)}\frac1{\varepsilon^{n+1}}
	\end{equation*}
The difference $F^{[\le \infty]}-F^{[\le H]}=F^{[>H]}$ is the generating function where the height is \emph{larger} than $H$. This needs to be normalized, in the form
\begin{equation*}
\frac{F^{[\le \infty]}-F^{[\le H]}}{F^{[\le \infty]}}\sim1-\Bigl(1-\frac{q}{\sqrt{p(p+4q)}}(p\varepsilon)^{H+1}\Bigr)^n
\end{equation*}
The generating function
\begin{equation*}
	\sum_{H\ge0}F^{[>H]}\sim	\sum_{H\ge0}\bigg[1-\Bigl(1-\frac{q}{\sqrt{p(p+4q)}}(p\varepsilon)^{H+1}\Bigr)^n\biggr]
\end{equation*}
is relevant for the \emph{average height}. The exponential approximation $(1-a)^n\approx \exp(-a\cdot n)$ for small $a$ plays a role here as well.
We use the following approximations
\begin{gather*}
p\varepsilon\sim\sqrt{p}+\frac p2+\dots,\ 
\frac{q}{\sqrt{p(p+4q)}}\sim\frac{\sqrt{p}}{2}+\frac{3}{16}p^{3/2}+\dots, \\ 
\frac{q}{\sqrt{p(p+4q)}}(p\varepsilon)^{H+1}\sim \frac{\sqrt{p}}{2}\sqrt{p}^{H+1}=\frac{p^{\frac{H+2}{2}}}{2}
\end{gather*}
and get
\begin{equation*}
\sum_{H\ge0}\bigg[1-\Bigl(1-\frac{p^{\frac{H+2}{2}}}{2}\Bigr)^n\biggr]\sim
\sum_{H\ge0}\bigg[1-\exp{\Bigl(-\frac{np^{1+{H}/{2}}}{2}\Bigr)}\biggr].
\end{equation*}
Set $N:=pn/2$, $\omega:=\sqrt{p}$ and rewrite:
\begin{equation*}
	\sum_{H\ge1}\Big[1-\exp{\bigl(-N\omega^{H}\bigr)}\Bigr].
\end{equation*}
From the enumeration by the Mellin transformation method \cite{FGD} we know that
\begin{equation*}
1-e^{-x}=\frac1{2\pi i}\int_{-\frac12-i\infty}^{-\frac12+i\infty}\Gamma(s)x^{-s}ds
\end{equation*}
and then
\begin{equation*}
	\sum_{H\ge1}[1-e^{-N\omega^H}]=\frac1{2\pi i}\int_{-\frac12-i\infty}^{-\frac12+i\infty}\Gamma(s)N^{-s}\sum_{H\ge1}\omega^{-Hs} ds=\frac1{2\pi i}\int_{-\frac12-i\infty}^{-\frac12+i\infty}\Gamma(s)N^{-s}\frac{1}{\omega^s-1} ds
\end{equation*}
The residue of the integrand at $s=0$ is of interest: $-\frac{\log N}{\log \omega}-\frac12-\frac{\gamma}{\log(\omega)}$.
Pulling this residue out  leads to the main asymptotic term of our sum, namely
\begin{equation*}
\frac{\log N}{\log \omega}=-\frac{\log \frac{pn}{2}}{\log \sqrt{p}}\sim -2\log_pn.
\end{equation*}
Lower order terms involve periodic fluctuations, see the discussion in \cite{Moran, GP97}.

\bibliographystyle{plain}

\end{document}